\documentclass[11pt]{article}
\usepackage{algorithm}
\usepackage{algorithmic}
\usepackage[left=2.5cm,top=1.5cm,right=2.5cm,nohead]{geometry}
\usepackage{graphicx}
\usepackage{multirow}

\title{An FMM Based on Dual Tree Traversal for Many-core Architectures}
\author{\textbf{R. Yokota}\\ King Abdullah University of Science and Technology}
\date{}
\begin{document}
\maketitle

\begin{abstract}
The present work attempts to integrate the independent efforts in the fast N-body community to create the fastest N-body library for many-core and heterogenous architectures. Focus is placed on low accuracy optimizations, in response to the recent interest to use FMM as a preconditioner for sparse linear solvers. A direct comparison with other state-of-the-art fast $N$-body codes demonstrates that orders of magnitude increase in performance can be achieved by careful selection of the optimal algorithm and low-level optimization of the code. The current N-body solver uses a fast multipole method with an efficient strategy for finding the list of cell-cell interactions by a dual tree traversal. A task-based threading model is used to maximize thread-level parallelism and intra-node load-balancing. In order to extract the full potential of the SIMD units on the latest CPUs, the inner kernels are optimized using AVX instructions.
\end{abstract}

\section{Introduction}
N-body simulations have traditionally been used in fields such as astrophysics and molecular dynamics, where the physics is naturally described by pairwise interaction of discrete bodies. However, applications for N-body solvers can be found in many other areas of physics including elastics, fluid dynamics, electromagnetics, acoustics, and quantum mechanics. This is made possible by transforming the partial differential equation governing the continuum field into an integral equation over discrete quadrature points. Therefore, N-body solvers are used in numerous scientific application codes and are considered as one of the key algorithmic components in scientific computing.

The fast multipole method (FMM) is a fast algorithm that reduces the complexity of the N-body problem from $\mathcal{O}(N^2)$ to $\mathcal{O}(N)$. It is widely regarded as one of the top 10 algorithms in scientific computing \cite{Board2000}, along with FFT and Krylov subspace methods. It is quite common that algorithms with low Byte/flop (dense linear algebra) have high complexity $\mathcal{O}(N^3)$, and algorithms with low complexity (FFT, sparse linear algebra) have high Byte/flop. The FMM has an exceptional combination of $\mathcal{O}(N)$ complexity and a Byte/flop that is even lower than matrix-matrix multiplication \cite{Yokota2012}. In other words, it is an \textit{efficient} algorithm that is \textit{compute} bound, which makes it an interesting alternative algorithm for many elliptic PDE solvers, on architectures of the future that will most likely have low Byte/flop.

The recent trend in computer architecture is shifting towards less and less Byte/flop ratios. A summary of the Byte/flop of flagship architectures from each vendor (as of August 2012) is shown in Table \ref{tab:byteperflop}. From this table, one can see that the Byte/flop will soon fall below 0.2, which makes it very difficult for most algorithms to extract a high percentage of peak flop/s from these architectures. According to the roofline study by Williams \textit{et al.} \cite{Williams2009}, sparse matrix-vector multiplication has about $4-5.88$ Byte/flop, stencil calculations have about $2-3$ Byte/flop, and 3-D FFT has about $0.61-0.92$ Byte/flop. It is clear how much the modern architectures are \textit{off balance} compared to the requirements of these common algorithms. This trend will most likely continue, so it is necessary to rethink the underlying algorithms in scientific simulations. Matrix-free methods seem to have an advantage in this respect, and FMM can be viewed as one of them.

In a recent study we have shown that FMM becomes faster than FFT when scaling to thousands of GPUs \cite{Yokota2012c}. The comparative efficiency between FMM and FFT can be explained from the asymptotic amount of communication. On a distributed memory system with $P$ nodes, a 3-D FFT requires two global transpose communications between $\sqrt{P}$ processes so the communication complexity is $\mathcal{O}(\sqrt{P})$. On the other hand, the hierarchical nature of the FMM reduces the amount of communication to $\mathcal{O}(\log P)$. A preliminary feasibility study for Exascale machines \cite{Gahvari2010} indicates that the necessary bandwidth for FFT could only be provided by a fat-tree or hypercube network. However, constructing such network topologies for millions of nodes is prohibitive in terms of cost, and the current trend of using torus networks is likely to continue. Therefore, network topology is another area where the trend in hardware is deviating from the requirements of common algorithms. Hierarchical methods are promising in this respect, and FMM is undoubtedly one of them.

It is clear from the above arguments that FMM could be an efficient alternative algorithm for many scientific applications on Exascale machines. One common objection is that FMM requires much more operations than other fast algorithms like multigrid and FFT, and therefore is much slower. However, as future microarchitectures move towards less and less Byte/flop, the asymptotic constant of the \textit{arithmetic} complexity becomes less of a concern. Therefore, the advantage in the \textit{communication} complexity trumps the disadvantages as mentioned in the previous paragraph.

\begin{table}[b]
\caption{Byte/flop of modern microprocessors. Byte/s is the theoretical peak of the memory bandwidth, and flop/s is the theoretical peak of the double-precision arithmetic throughput when fused-multiply-add units are fully utilized, along with maximum clock frequency by turbo boost (if available).}
\label{tab:byteperflop}
\begin{center}
\begin{tabular}{|c|c|c|c|c|c|}
\hline
Vendor & Microarchitecture & Model & Byte/s & Flop/s & Byte/flop\\
\hline \hline
Intel & Sandy Bridge & Xeon E5-2690 & 51.2 & 243.2 & 0.211\\
\hline
AMD & Bulldozer & Opteron 6284 SE & 51.2 & 217.6 & 0.235\\
\hline
AMD & Southern Islands & Radeon HD7970 (GHz Ed.) & 288 & 1010 & 0.285\\
\hline
NVIDIA & Fermi GF110 & Tesla M2090 & 177 & 665 & 0.266\\
\hline
IBM & PowerPC & PowerPC A2 (BG/Q) & 42.6 & 204.8 & 0.208\\
\hline
Fujitsu & SPARC64 & SPARC64 IXfx (FX10) & 85 &  236.5 & 0.359\\
\hline
\end{tabular}
\end{center}
\end{table}

\section{Related Work}
FMM is a relatively new algorithm compared to well established linear algebra solvers and FFT, and has ample room for both mathematical and algorithmic improvement, not to mention the need to  develop highly optimized libraries. We will briefly summarize the recent efforts in this area by categorizing them into; CPU optimization, GPU optimization, MPI parallelization, algorithmic comparison, auto-tuning, and data-driven execution.

In terms of CPU optimization, Chandramowlishwaran \textit{et al.} \cite{Chandramowlishwaran2010} have done extensive work on thread-level parallelization with NUMA-aware optimizations and manual SIMD vectorization. They performed benchmarks on four different microarchitectures; Harpertown, Barcelona, Nehalem-EP, and Nehalem-EX, and exploit the simultaneous multithreading features of the Nehalem processors. With the various optimizations, they were able to double the performance on Harpertown and Barcelona processors, while obtaining a 3.7x speed-up on Nehalem-EX, compared to their previous best multithreaded and vectorized implementation. Some of the optimizations are specific to the kernel-independent FMM (KIFMM) \cite{Ying2004}, and do not apply to FMMs based on Cartesian, spherical harmonics, or planewave expansions. It would be interesting to compare highly optimized kernels for these different types of expansions (including KIFMM) to determine at what order of expansion they cross over.

For GPU optimization for uniform distributions, Takahashi \textit{et al.} \cite{Takahashi2012} developed a technique to turn the originally memory bound M2L kernel into a compute bound GPU kernel by exploiting the symmetry of the interaction stencil. There are alternative techniques to make the M2L kernel compute bound on GPUs, by recalculating the M2L translation matrix on-the-fly \cite{Yokota2011c}. However, the fact that the technique by Takahashi \textit{et al.} precalculates the M2L translation matrix and still allows the M2L kernel to be compute bound, means it is that much more efficient. Increasing the number of particles per leaf \cite{Yokota2009} or calculating the M2L kernel on the CPU \cite{Hu2011} is an easy optimization strategy for heterogenous architectures. The value of this work lies in the fact that they tackle the GPU-M2L optimization problem head on, and actually succeed in making it compute bound. One limitation of their method is that the particle distribution must be somewhat uniform in order to exploit the symmetry of the M2L translation stencil. Though, classical molecular dynamics is an important application that \textit{will} be able to benefit from this method.

For GPU optimization of non-uniform distributions for low accuracy, Bedorf \textit{et al.} \cite{Bedorf2012} developed a highly optimized treecode \texttt{bonsai}\footnote{https://github.com/treecode/Bonsai} that runs entirely on GPUs. For three significant digits of accuracy in the force calculation for the Laplace kernel, their code is an order of magnitude faster than any other treecode or FMM on GPUs. Their code takes less than 50 milliseconds to calculate one million particles (including the tree construction) on an NVIDIA Fermi GF110 series GPU for the aforementioned accuracy, whereas other codes take about 500 milliseconds to achieve the same accuracy for one million particles. By using warp unit execution and \texttt{\_\_ballot()} functions they have eliminated all synchronization barriers from their CUDA kernel. Furthermore, by using mask operations and prefix sums, they have eliminated all conditional branching from their code. However, an application which requires high accuracy may find this code difficult to use, since the order of expansion is not adjustable past $p>3$. Nonetheless, anyone who has experience optimizing fast N-body codes knows that optimizing for low accuracy is much more challenging than optimizing high accuracy codes, so the impact of this work is remarkable.

For MPI parallelization of an adaptive tree, Winkel \textit{et al.} \cite{Winkel2012} propose a request based point-to-point communication scheme that dedicates a thread for communication on each node, while the evaluation is done using the remaining threads simultaneously, in their code \texttt{PEPC}\footnote{https://trac.version.fz-juelich.de/pepc}. Meanwhile, Jetley \textit{et al.} \cite{Jetley2010} used \texttt{CHARMM++} to dynamically load-balance on large scale heterogenous environments, obtaining large gains by task aggregation and memory pooling, in their code \texttt{ChaNGa}\footnote{http://software.astro.washington.edu/nchilada}. These efforts to provide an alternative to bulk-synchronous communication models could eventually become advantageous, but the lack of comparison against current state-of-the-art bulk-synchronous methods \cite{Rahimian2010,Yokota2012c} that also scale to the full node of the largest supercomputers, makes it difficult to predict exactly when they will become advantageous (if they do at all).

As another aspect of MPI parallelization of an adaptive tree, Winkel \textit{et al.} use a 64-bit global Morton key and partition the domain to equally distribute the number of keys, while using the traversal time of the previous step as weights. This technique will only allow a maximum depth of the global tree of 21 levels before the 64-bit key overflows. Current large scale simulations already exceed this depth \cite{Rahimian2010}, and the method of Winkel \textit{et al.} will require the use of multiple integers to store the key in such cases. Looking into the future, it would seem that the use of ever more integers to store the key will slow down the sorting of these keys, but the fact that Bedorf \textit{et al.} \cite{Bedorf2012} can sort over a billion 96-bit keys per second using the \texttt{back40computing} library\footnote{http://code.google.com/p/back40computing} means that this remains a viable option. It would be interesting to compare this approach with alternative techniques that eliminate the use of global keys altogether \cite{Yokota2012c}.

For algorithmic comparison, Fortin \textit{et al.} \cite{Fortin2011} compare the BLAS based adaptive FMM by Coulaud \textit{et al.} \cite{Coulaud2008},  \texttt{NEMOtree1}\footnote{http://bima.astro.umd.edu/nemo} implementation of the Barnes \& Hut treecode \cite{Barnes1986}, the \texttt{GADGET-2}\footnote{http://www.mpa-garching.mpg.de/gadget} code by Springel \textit{et al.} \cite{Springel2005} (without using the treePM part), and \texttt{falcON}\footnote{http://bima.astro.umd.edu/nemo/man\_html/gyrfalcON.1.html} code by Dehnen \cite{Dehnen2002}. The comparison is for 6 different datasets with various non-uniformity, and the time-to-solution and memory consumption was compared for different problem sizes and various number of cores, while the accuracy was set to roughly 3 significant digits for the force. There was over an order of magnitude difference in the speed between these codes, \texttt{falcON} being consistently the fastest for all types of non-uniform distributions, and exhibiting perfect $\mathcal{O}(N)$ complexity throughout the examined range of $N$. A major drawback of \texttt{falcON} is that it is not parallel. In the present work, we extended the \texttt{falcON} framework to MPI, pthreads, OpenMP, Intel TBB, SSE, and AVX based parallelism, as shown later in this article.

An example of auto-tuning efforts in FMMs is given by Yokota \& Barba \cite{Yokota2012}, where the dual tree traversal approach by Dehnen \cite{Dehnen2002} was combined with an auto-tuning mechanism by selecting between M2L and M2P kernels to produce a treecode-FMM hybrid code, and by selecting between M2L and P2P kernels to optimize the number of particles at the leaf at which P2P is performed. This is an advantage when developing a black-box library that runs on both CPUs and GPUs without having to worry about selecting the optimal number of particles per leaf for each architecture. This code -- \texttt{exaFMM}\footnote{https://bitbucket.org/rioyokota/exafmm-dev} -- inherits all the techniques of \texttt{falcON} such as the dual tree traversal, mutual interaction, error aware local optimization of the opening angle $\theta$, and a multipole acceptance criterion based on $B_{max}$ and $R_{max}$ \cite{Dehnen2002}. At the same time, it adds features such as arbitrary order of expansion by using template metaprogramming, periodic boundary conditions, CUDA version of all FMM kernels, and a scalable MPI implementation. The present work is a further extension of this code to thread-level and SIMD parallelism.

Last but not least, there have been a few recent attempts to use event-driven runtime systems for optimizing the data flow of parallel FMM codes. Ltaief \& Yokota \cite{Ltaief2012} used \texttt{QUARK}\footnote{http://icl.cs.utk.edu/quark/index.html} to schedule the threads dynamically according to a directed acyclic graph that represents the data flow of \texttt{exaFMM}'s dual tree traversal. Dekate \textit{et al.} take a similar approach but using \texttt{ParalleX}\footnote{http://stellar.cct.lsu.edu/tag/parallex} for the runtime system and a classical Barnes-Hut treecode \cite{Barnes1986} for the fast $N$-body solver. Agullo \textit{et al.} use \texttt{StarPU}\footnote{http://runtime.bordeaux.inria.fr/StarPU} along with the black-box FMM \cite{Fong2009} on a CPU+GPU environment. Peric\'as \textit{et al.} use \texttt{OmpSs} \footnote{http://pm.bsc.es/ompss} with \texttt{exaFMM}. These are all preliminary results, and large gains from using runtime systems have yet to be reported.

In summary, the present work aims to provide the best implementation of the best algorithm for fast $N$-body solvers on many-core and heterogeneous architectures based on the thorough investigation of existing algorithms and implementations mentioned above. Our results indicate that many orders of magnitude in acceleration can be achieved by integrating many of these techniques that have not been previously combined. Upon integrating these various techniques, it was necessary to first reverse engineer many existing codes and to systematically investigate the individual contributions of each algorithmic improvement and each optimization trick that was not mentioned in the accompanying publication. Then, the value of these individual contributions had to be reevaluated based on predictions of how they will preform on many-core co-processors and GPUs. Finally, careful design decisions had to be made on which techniques to adopt, and what standardized constructs could be used to implement them in the simplest and most portable way in the integrated code. This ongoing work is being made available as an open source FMM library for Exascale -- \texttt{exaFMM}, which is currently hosted on bitbucket\footnote{https://bitbucket.org/rioyokota/exafmm-dev}. Although the current work is influenced by many other predecessors, the work cited in this section are the ones that influenced the design decision of the final integrated code the most.

\section{Present Method and Results}
There are many variants of fast $N$-body methods, as described in the previous section. The argument in the previous section is strengthened by providing a more detailed description of the underlying techniques. The significance of the present work can be understood in more depth if these subtle differences in the algorithms and their implications on future architectures are explained clearly. Each subsection within this section revisits certain aspects of fast $N$-body methods in light of many-core and heterogenous architectures. In the end, we will compare 5 major implementations of treecodes and FMMs to show the relative performance of the resulting code.

\begin{table}[b!]
\caption{Asymptotic complexity of the storage and arithmetic with respect to the order of expansion $p$ for the different expansions in fast $N$-body methods (3-D). There are two references for each expansion, the first is the original work describing the use of each expansion in the context of FMM, while the second has the simplest mathematical description of the translation operations.}
\label{tab:expansion}
\begin{center}
\begin{tabular}{|c|c|c|}
\hline
Type of expansion (+M2L acceleration) & Storage & Arithmetic\\
\hline \hline
Cartesian Taylor \cite{Applequist1989,Visscher2010} & $\mathcal{O}(p^3)$ & $\mathcal{O}(p^6)$\\
\hline
Cartesian Chebychev \cite{Dutt1996,Fong2009} & $\mathcal{O}(p^3)$ & $\mathcal{O}(p^6)$\\
\hline
Spherial harmonics \cite{Greengard1988,Cheng1999} & $\mathcal{O}(p^2)$ & $\mathcal{O}(p^4)$\\
\hline
Spherial harmonics+rotation \cite{White1996,Cheng1999} & $\mathcal{O}(p^2)$ & $\mathcal{O}(p^3)$\\
\hline
Spherial harmonics+FFT \cite{Elliott1996,Kurzak2006} & $\mathcal{O}(p^2)$ & $\mathcal{O}(p^2\log^2 p)$\\
\hline
Planewave \cite{Greengard1997,Kurzak2006} & $\mathcal{O}(p^2)$ & $\mathcal{O}(p^3)$\\
\hline
Equivalent charges \cite{Anderson1992,Makino1999} & $\mathcal{O}(p^2)$ & $\mathcal{O}(p^4)$\\
\hline
Equivalent charges+FFT \cite{Ying2004, Langston2011} & $\mathcal{O}(p^3)$ & $\mathcal{O}(p^3\log p)$\\
\hline
\end{tabular}
\end{center}
\end{table}

\subsection{Choice of Series Expansions}
The accuracy of fast $N$-body algorithms is controllable by adjusting the order of multipole/local expansions. There are many series expansions that can be used for approximating Green's functions, \textit{e.g.} Cartesian, spherical harmonics, planewaves, and equivalent charges. These series expansions have different characteristics in terms of Byte/flop, parallelism, and asymptotic complexity. A majority of the existing techniques to accelerate FMMs are based on old programming paradigms from back when arithmetic operations were the bottleneck. For example, pre-computation of translation matrices will reduce the arithmetic workload of \textit{computing} the matrix entries multiple times, but requires extra bandwidth for \textit{loading} the matrix every time. This goes against the current trend in hardware, where Byte/flop is decreasing at every level of the memory hierarchy. Furthermore, it may even be necessary to reconsider the choice of series expansions based on how regular the data access patterns are, and how parallel the operations are.

The asymptotic complexity of the storage and arithmetic operation for the M2L kernel for each expansion is shown in Table \ref{tab:expansion}. Cartesian Taylor expansions are most widely used in treecodes, while other types of expansions are only used in FMMs. The common (but misleading) claim that treecodes are faster than FMMs for low accuracy, is actually a consequence of the choice of series expansions. Dehnen \cite{Dehnen2002} used Cartesian expansions for FMM, and showed that it was an order of magnitude faster than the FMM using spherical harmonics or planewaves by Cheng \textit{et al.} \cite{Cheng1999} for low accuracy, but was slower for higher accuracy. This result clearly demonstrates that it is not the use of treecodes that makes the difference, but rather the Cartesian Taylor expansion that increases the relative performance for lower accuracy. The FMM is always faster than treecodes for \textit{any} accuracy if the only difference was cell-cell interactions vs. cell-particle interactions, and all other conditions (choice of expansion, particles per leaf cell, opening angle $\theta$) remain the same \cite{Yokota2012}.

Cartesian expansions have higher arithmetic complexity but smaller asymptotic constants, whereas methods with lower arithmetic complexity tend to have larger asymptotic constants. For example, White and Head-Gordon \cite{White1996} compare the spherical harmonics expansion with rotation $\mathcal{O}(p^3)$ and without rotation $\mathcal{O}(p^4)$, and observe only a 2 fold speed up even for $p=21$. When comparing the break-even point between different series expansions, it is necessary to consider the low order terms as well. For example a $\mathcal{O}(p^6)$ method actually has a complexity of $ap^6+bp^5+cp^4+dp^3+ep^2+fp+g$. If the coefficients of the lower order terms are large, a $\mathcal{O}(p^6)$ method may exhibit lower order asymptotic behavior for practical accuracy.

\begin{figure}[t]
\centering
\includegraphics[width=0.5\textwidth]{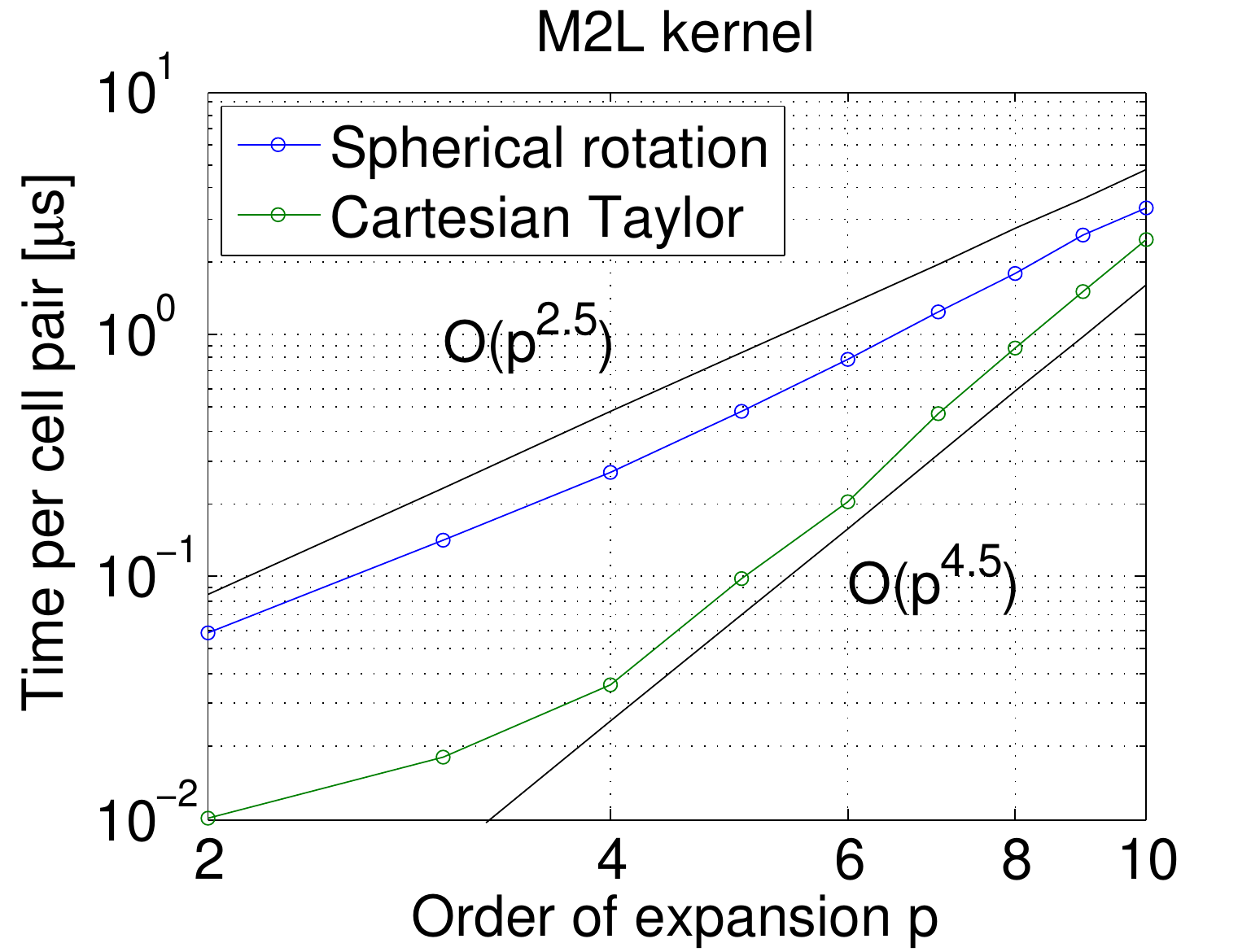}
\caption{Calculation time of M2L translation between a single pair of cells for different orders of expansion $p$. The time is shown in microseconds. This data is calculated by performing a full FMM calculation and timing the M2L kernel while counting the number of cell pairs that call the M2L kernel. Calculations are performed on a single socket of a Xeon X5650 (2.67GHz, 6 physical cores).}
\label{fig:orderp}
\end{figure}

In order to demonstrate this behavior, we plot the calculation time of a single M2L kernel for a pair of cells against the order of expansion $p$ for the $\mathcal{O}(p^3)$ rotation based spherical harmonics expansion and $\mathcal{O}(p^6)$ Cartesian Taylor expansion in Figure \ref{fig:orderp}. The time is shown in microseconds, and the runs were performed on a single socket of a Xeon X5650 (Westmere-EP, 2.67 GHz, 6 physical cores). The timings are measured using an actual FMM calculation with $N=1,000,000$ particles, while timing the M2L kernels and counting the number of M2L kernel calls. The timing per M2L kernel is obtained by dividing the total M2L time by the number of M2L kernel calls. Even though the leading order term of the complexity of these methods are $\mathcal{O}(p^3)$ and $\mathcal{O}(p^6)$, they are closer to $\mathcal{O}(p^{2.5})$ and $\mathcal{O}(p^{4.5})$ in the practical accuracy range. Furthermore, the present implementation of the rotation based spherical harmonics expansion seems to crossover with the Cartesian Taylor expansion at around $p\approx12$. Our Cartesian Taylor kernels are highly optimized using template meta-programming with specialization for low $p$. On the other hand, our rotation based spherical harmonics kernels are not optimized to the same extent, and simply use precalculated Wigner rotation matrices and translation matrices. Both versions use OpenMP but no SIMD intrinsics.

It would be useful to the FMM community to complete Figure \ref{fig:orderp} with the other entires in Table 2, and for a broader range of $p$. However, such results are not only hardware dependent, but also highly implementation dependent. The tests would have to be run on a wide range of architectures, while tuning each kernel for each architecture. Furthermore, low level optimization of kernels of this complexity require a significant effort, and the absence of highly tuned open source implementations of these kernels makes it difficult to perform an exhaustive comparison between different expansions. Therefore, completing Figure \ref{fig:orderp} with the remaining series expansions is an ongoing investigation, which we hope to report in a future publication.

In the current work, we will use Cartesian Taylor expansions and focus on low accuracy FMMs. High accuracy FMMs with larger order of expansions will yield more fine grained parallelism, and will naturally offer better performance on many-core CPUs and GPUs. The real challenge is to optimize low accuracy FMMs on these architectures, and the reward will be substantial if this is possible. Recent efforts to view the FMM as a general elliptic PDE solver \cite{Langston2011} have opened the possibility to use it as a preconditioner for even a broader range of applications. If the FMM is to be used as a preconditioner the required accuracy would not be high, which is the reason why we are interested in low accuracy FMM on future architectures. In terms of relative performance compared to multigrid preconditioners, the treecode by Bedorf \textit{et al.} \cite{Bedorf2012} can solve 8 levels in roughly 50 milliseconds on 1 GPU \footnote{The latest code from https://github.com/treecode/Bonsai is much faster than the one reported in their paper}, while the multigrid code by Goddeke \textit{et al.} also takes about 50 milliseconds for 8 levels on 1 GPU \cite{Goddeke2011} to achieve the same accuracy.

\subsection{Treecodes and FMMs}
The two variants of hierarchical $N$-body methods -- treecodes and FMMs -- have followed very different paths of evolution over the past two decades. The primary reason for this divergence is the difference in the communities that adopted them. Treecodes have been used predominantly in the astrophysics community, where the distribution of stars are highly non-uniform, variation of local time scales is large, and the required accuracy is relatively low. On the other hand, FMMs received more attention from the applied mathematics community, and more emphasis was placed on asymptotically faster techniques with respect to the order of expansion $p$, extension to other kernels besides Laplace, and other generalizations to make them applicable to a wider range of problems. The variety of series expansions mentioned in the previous subsection is one example where FMMs are more advanced. Conversely, an area where treecodes are more sophisticated is the handling of adaptive tree structures for highly non-uniform distributions, which was an absolute requirement for their primary application. There are many interesting design decisions that arise when merging these orthogonal efforts into a unified framework.

The key differences between typical treecodes and FMMs are summarized in Table \ref{tab:difference}. The first item is the difference in how the particles are grouped into a tree structure. Treecodes commonly use rectangular grouping which squeeze the bounds of each cell\footnote{In our descriptions, the terminology ``nodes" and ``cells" essentially mean the same thing, where ``node" refers to the elements in the tree data structure (from a computer science perspective), while ``cell" refers to the corresponding geometric object in physical space (from a scientific computing perspective).} to exactly fit the residing particles \cite{Duan2001}, while FMMs always use cubic cells. When combined with the second item in the table for a sophisticated definition of ``well-separated" cells, the use of tight rectangular cells results in an efficient error bounding mechanism. On the other hand, the use of cubic cells assumes a worst case scenario where the farthest particle from the cell center is at the corner of the cell. This may be a good approximation at coarser levels of a uniform tree, but is an extremely crude assumption at the leaf level of an adaptive tree. This assumption, along with the rigid definition of well-separated cells to be simply the non-neighboring cubes, results in an inefficient mechanism to achieve a given accuracy. This effect is much more prominent for low accuracy calculations, where the number of particles per leaf cell becomes small. 

\begin{table}[b]
\caption{Key differences between typical treecodes and FMMs.}
\label{tab:difference}
\begin{center}
\begin{tabular}{|c|c|c|}
\hline
& Treecode & FMM\\
\hline \hline
Octree node shape & rectangle & cube \\
\hline
Definition of  ``well-separated" & multipole acceptance criterion & non-neighboring cubes\\
\hline
Center of expansion & center of mass & geometric center of cell\\
\hline
Accuracy control & opening angle $\theta$ & order of expansion $p$\\
\hline
Data structure of tree & linked lists & Morton key\\
\hline 
\multirow{2}{*}{Finding interacting nodes} & \multirow{2}{*}{tree traversal} & parent's neighbor's children\\
& & that are non-neighbors \\
\hline
Domain decomposition & orthogonal recursive bisection & Morton key partitioning\\
\hline
\end{tabular}
\end{center}
\end{table}

\begin{figure}[t]
\centering
\includegraphics[width=1.0\textwidth]{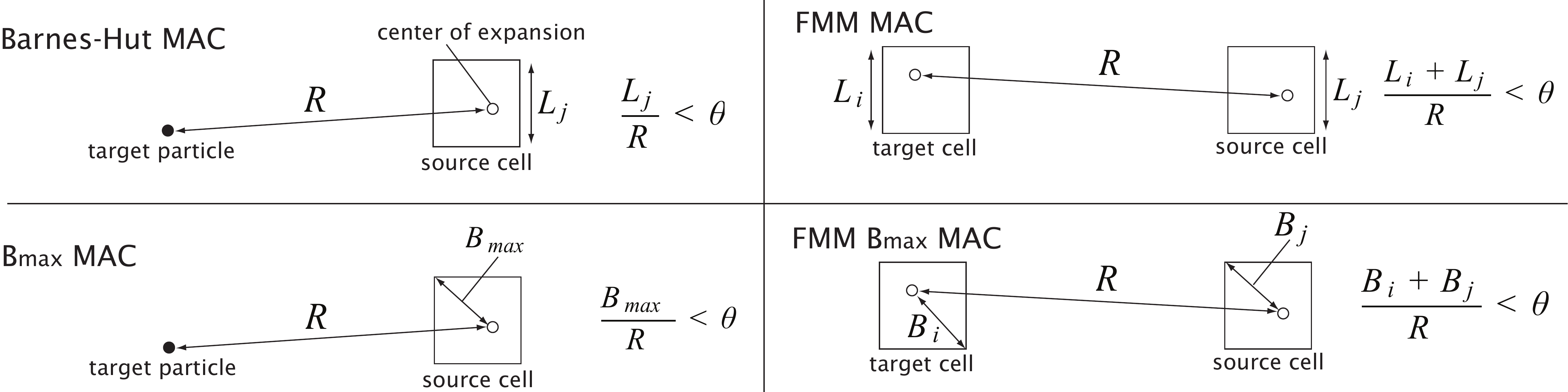}
\caption{Illustration of the different types of multipole acceptance criteria (MAC). Barnes-Hut MAC is the original one proposed by Barnes and Hut \cite{Barnes1986}, $B_{max}$ MAC is a tighter MAC proposed by Salmon and Warren \cite{Salmon1994a}, and FMM MAC is the extension of the $B_{max}$ MAC to FMM \cite{Warren1995}.}
\label{fig:mac}
\end{figure}

The multipole acceptance criterion (MAC) for treecodes has many variants as described by Salmon and Warren \cite{Salmon1994a}. The concept of MACs have also been extended to FMM by the same authors \cite{Warren1995} where they provide error bounds for the MAC based FMM. The different types of MAC are shown in Figure \ref{fig:mac}. The Barnes-Hut MAC is the original one proposed by Barnes and Hut \cite{Barnes1986}. It sets the threshold for accepting the multipole expansion according to the ratio between the cell size and distance between the target particle and center of expansion. Cells that are not accepted are further divided into their child cells and reevaluated. Therefore, the accuracy of treecodes can be controlled by adjusting the opening angle $\theta$. Though, decreasing $\theta$ too much will ultimately turn the treecode into a $\mathcal{O}(N^2)$ direct summation method (with no error by definition). The $B_{max}$ MAC gives a tighter error bound \cite{Salmon1994a} and the same concept can be extended to FMMs \cite{Warren1995}.

Using the notation in Figure \ref{fig:mac}, the leading error term in the series expansion for the M2L kernel is of the order $\mathcal{O}(\left(\frac{B_i+B_j}{R}\right)^p)$. Therefore, the FMM has a simple error bound of $\mathcal{O}(\theta^p)$ \cite{Dehnen2002}. This error bound is easily extendable to treecodes by assuming $B_i=0$, which is equivalent to assuming that the target cell is just one particle and has no size. As the third item in Table \ref{tab:difference} indicates, treecodes increase the accuracy by decreasing $\theta$, while keeping $p$ constant ($p=3$ in most cases). Conversly, FMMs increase the accuracy by increasing $p$, while keeping $\theta$ constant (most commonly $\theta=1$ according to the definition of ``FMM MAC" in Figure \ref{fig:mac}). At every level of the tree, the M2L interaction list forms a constant-thickness layer on the surface of a sphere with radius $\frac{B_i+B_j}{\theta}<R<2\frac{B_i+B_j}{\theta}$, which means the lists size grows as $\mathcal{O}(\theta^{-2})$ (surface of sphere with normalized radius $1/\theta$). If one were to use a $\mathcal{O}(p^3)$ method for the M2L kernel, the complexity with respect to $p$ and $\theta$ will be $\mathcal{O}(\theta^{-2}p^3)$. There should exist an optimal balance between $p$ and $\theta$ that minimizes the $\mathcal{O}(\theta^{-2}p^3)$ work, while keeping the $\mathcal{O}(\theta^p)$ error constant.

\begin{table}[b!]
\caption{Balance between order of expansion $p$ and opening angle $\theta$ for Cartesian Taylor expansions to achieve a specific error tolerance. $Err$ is the relative $L^2$ norm error of the force (not potential) for the Laplace kernel, where the result from a direct summation is used as the reference value. The time is measured for the dual tree traversal (M2L kernels+P2P kernels+list construction) for $N=10^5$ particles for one time step. The maximum number of particles per leaf cell is set to $N_{crit}=30$ for all cases.}
\label{tab:ptheta}
\begin{center}
\begin{tabular}{|c|c|c|c|c|}
\hline
& $p=3$ & $p=4$ & $p=5$ & $p=6$\\
\hline \hline
\multirow{2}{*}{$Err=10^{-2}$} &
$\theta=1.00$ & $\theta=1.18$ & $\theta=1.23$ & $\theta=1.24$\\
& $time=0.016s$ & $time=0.012s$ & $time=0.015s$ & $time=0.026s$\\
\hline
\multirow{2}{*}{$Err=10^{-3}$} &
$\theta=0.67$ & $\theta=0.78$ & $\theta=0.91$ & $\theta=0.94$\\
& $time=0.036s$ & $time=0.027s$ & $time=0.024s$ & $time=0.038s$\\
\hline
\multirow{2}{*}{$Err=10^{-4}$} &
$\theta=0.30$ & $\theta=0.49$ & $\theta=0.62$ & $\theta=0.70$\\
& $time=0.22s$ & $time=0.085s$ & $time=0.071s$ & $time=0.073s$\\
\hline
\multirow{2}{*}{$Err=10^{-5}$} &
$\theta=0.12$ & $\theta=0.20$ & $\theta=0.36$ & $\theta=0.45$\\
& $time=1.38s$ & $time=0.59s$ & $time=0.21s$ & $time=0.21s$\\
\hline
\end{tabular}
\end{center}
\end{table}

The various combinations of $p$ and $\theta$ for a given accuracy are shown in Table \ref{tab:ptheta}. This table was generated by first setting a target accuracy (relative $L^2$ norm error of the force (not potential) for the Laplace kernel, shown as $Err$) and setting $p$. Then, we searched through the parameter space of $\theta$ to find the largest value that achieves the target accuracy. These, values of $\theta$ are shown in the table along with the amount of time it takes to perform the dual tree traversal for each case. The tree construction, P2M, M2M, L2L, L2P kernels are not included in these timings. The number of particles is $N=10^5$, and the maximum number of particles per leaf cell is set to $N_{crit}=30$ for all cases. Here, the definition of $\theta$ is a more advanced version than the ones shown in Figure \ref{fig:mac}, and uses both $B_{max}$ and $R_{max}$ \cite{Dehnen2002} but without error optimization of $R_{max}$. The optimal $p$ shifts to larger values when more accuracy is required. According to the current definition, $Err=10^{-5}$ means that 5 significant digits match with the results of a direct summation, which is sufficient for many applications in fluid dynamics, structural dynamics, and molecular dynamics. By comparing with Figure \ref{fig:orderp}, one can conclude that Cartesian Taylor expansions are much faster than spherical harmonics with rotation for this range of $p$.

\subsection{Dual tree traversal}
As summarized in Table \ref{tab:difference}, the underlying data structure in treecodes and FMMs is quite different. Treecodes usually construct a linked-list, where each node is expressed as a structure which contains a pointer to it's children. By looping over it's children in a recursive function call (or by using a stack or queue), all nodes can be visited once in a topdown manner, \textit{i.e.} a tree traversal. The multipoles in the well-separated cells are easily found by applying the MAC during the tree traversal. This procedure is shown in Algorithm \ref{al:evaluateTreecode}.

On the other hand, most FMM codes do not construct linked-lists between the tree nodes, nor do they traverse a tree. Instead, they bin the particles based on Morton/Hilbert keys, and loop over all the target cells in the tree (level wise), and explicitly construct a list of source cells that are well-separated. This procedure does not require a linked tree structure. The definition of well-separated cells is usually the ``target cell's parent's neighbor's children that are not direct neighbors of the target cell". Parents and children are found by multiplying and dividing the Morton key by 8. Neighbors are found by deinterleaving the bits in the key into a three-dimensional index, incrementing/decrementing in a certain direction, and interleaving them back to a single key. In adaptive trees these keys are mapped to a unique index by hashing (or simply arrays with zeros for the empty boxes), and this unique index is used to load/store the corresponding multipole/local expansion coefficients. This procedure is shown in Algorithm \ref{al:evaluateFMM}.

\begin{table}[b!]
\caption{Different types of tree traversals and interaction list generation.}
\label{tab:traversal}
\begin{center}
\begin{tabular}{|c|c|c|c|c|}
\hline
& Treecode & FMM & FMM (adaptive) & FMM (DTT)\\
\hline \hline
Data structure & linked list & array & array & linked list\\
\hline
Complexity & $\mathcal{O}(N\log N)$ & $\mathcal{O}(N)$ & $\mathcal{O}(N)$ & $\mathcal{O}(N)$\\
\hline
Interaction list & implicit & explicit & explicit & implicit\\
\hline
Cell shape & rectangular & cubic & cubic & rectangular\\
\hline
$\theta$ adjustable & yes & no & no & yes\\
\hline
Parallelism & for loop & for loop & for loop & task-based\\
\hline
\end{tabular}
\end{center}
\end{table}

There are more advanced techniques for finding the interaction lists in FMMs. Cheng \textit{et al.} extended the originally non-adaptive interaction list with only 2 categories ``well-separated" and ``neighbor", to lists for adaptive trees using 4 categories \cite{Cheng1999}. Gumerov and Duraiswami proposed a way to reduce the number of interactions for the ``well-separated" cells of the original non-adaptive interaction list from 189 to 119 without sacrificing the accuracy \cite{Gumerov2008}. However, both of these techniques have the same structure of using the ``parent's neighbors" to define the well-separated cells, and finding neighbors in this way have several disadvantages. First, special treatment of cells on the boundaries is required since they will not have a full list of 26 neighbors. Second, the definition of ``well-separatedness" cannot be flexibly controlled, which makes a large difference as we have demonstrated in Table \ref{tab:ptheta}. Third, only cubic cells can be used because the cell's dimensions are used as units to measure well-separatedness, whereas rectangular cells have been preferred in treecodes for highly non-uniform distributions. This third disadvantage also has implications for load-balancing adaptive global trees, which is an important issue that is beyond the scope of the present work.

It is possible to overcome the aforementioned disadvantages by using a dual tree traversal \cite{Warren1995, Dehnen2002}, where the target tree and source tree are traversed simultaneously to efficiently create the interaction list. A pseudocode for the dual tree traversal is shown in Algorithms \ref{al:evaluateFMM2} and \ref{al:interact}. Table \ref{tab:traversal} summarizes the characteristics of four different ways to systematically evaluate the far field and near field contributions in fast $N$-body methods. The dual tree traversal inherits the flexibility and adaptiveness of the tree traversal used in treecodes, while achieving $\mathcal{O}(N)$ complexity. The ``well-separatedness" can be defined flexibly and the use of rectangular cells are allowed. The target tree and source tree can be completely independent tree structures, since the distance between cells is identified from the coordinates of the center of expansion at each cell. The use of Morton/Hilbert keys is not a requirement, but can be used to accelerate the tree construction process. One drawback of the dual tree traversal could be the loss of obvious parallelism at the outermost loop, as can be seen by comparing Algorithms \ref{al:evaluateTreecode}, \ref{al:evaluateFMM}, and \ref{al:evaluateFMM2}. Treecodes have an outer loop over all target leafs, while conventional FMMs have an outer loop over all target cells. However, the dual tree traversal does not have any ``for loop" structures that could be readily parallelized by, for example a \texttt{\#pragma omp parallel for} directive. This may or may not be an issue since standards such as OpenMP 3.0 and Intel Thread Building Blocks (TBB) now provide task-based threading models, which do not require ``for loops" for parallelization. In this model, tasks are spawned by the user and dispatched to idle threads dynamically during runtime. Extracting thread level parallelism from the dual tree traversal on NUMA architectures is an interesting challenge for data-driven execution models such as \texttt{QUARK}, \texttt{StarPU}, \texttt{OmpSs}, and \texttt{Parallex}.

The lack of ``for loop" type parallelism in dual tree traversals makes it difficult to simply assign a certain part of the work to a thread. Therefore, on GPUs and many-core co-processors, it is crucial to use task-based parallelization tools available in CUDA 5.0, OpenMP 3.0 and Intel TBB. However, the lack of ``for loop" type parallelism in dual tree traversals itself is not a limitation, since these task-based parallelization tools now exist on the latest compilers. It is important to stay updated on the advancement in compiler technology and to continuously consider all variants of FMMs that can take advantage of them.

\begin{algorithm}
\caption{EvaluateTreecode()}
\label{al:evaluateTreecode}
\begin{algorithmic}
\FOR{all target leafs}
\STATE push source root cell to stack
\WHILE{stack is not empty}
\STATE pop stack
\FOR{all children of cell}
\IF{child cell is a leaf}
\STATE call P2P kernel
\ELSE
\IF{child cell satisfies MAC}
\STATE call M2P kernel
\ELSE
\STATE push child cell to stack
\ENDIF
\ENDIF
\ENDFOR
\ENDWHILE
\ENDFOR
\end{algorithmic}
\end{algorithm}

\begin{algorithm}
\caption{EvaluateFMM()}
\label{al:evaluateFMM}
\begin{algorithmic}
\FOR{all target cells}
\STATE find parent cell by dividing Morton key
\FOR{all children of parent cell's neighbor}
\STATE call M2L kernel
\ENDFOR
\ENDFOR
\FOR{all target leafs}
\FOR{all neighbor's}
\STATE call P2P kernel
\ENDFOR
\ENDFOR
\end{algorithmic}
\end{algorithm}

\begin{algorithm}
\caption{EvaluateFMMUsingDualTreeTraversal()}
\label{al:evaluateFMM2}
\begin{algorithmic}
\STATE push pair of root cells ($A$,$B$) to stack
\WHILE{stack is not empty}
\STATE pop stack to get ($A$,$B$)
\IF{target cell is larger then source cell}
\FOR{all children $a$ of target cell $A$}
\STATE Interact($a$,$B$)
\ENDFOR
\ELSE
\FOR{all children $b$ of source cell $B$}
\STATE Interact($A$,$b$)
\ENDFOR
\ENDIF
\ENDWHILE
\end{algorithmic}
\end{algorithm}

\begin{algorithm}
\caption{Interact($A$,$B$)}
\label{al:interact}
\begin{algorithmic}
\IF{$A$ and $B$ are both leafs}
\STATE call P2P kernel
\ELSE
\IF{$A$ and $B$ satisfy MAC}
\STATE call M2L kernel
\ELSE
\STATE push pair ($A$,$B$) to stack
\ENDIF
\ENDIF
\end{algorithmic}
\end{algorithm}

\subsection{Code optimization}
A highly optimized code must take advantage of SIMD instructions, to extract the full potential of modern microprocessors. For the current purpose of optimizing FMM for low accuracy, single-precision is sufficient. For this case SSE and AVX can offer a maximum of $\times4$ and $\times8$ acceleration, respectively. Furthermore, $N$-body kernels involve a reciprocal square root operation, which can be accelerated by an intrinsic function \texttt{rsqrt\_ps()} in SSE/AVX. Furthermore, by reordering the instructions to interleave the load operations with fused-multiply-add operations, one can achieve flop/s that is very close to the theoretical peak, given the low Byte/flop ratio of these $N$-body kernels.

\begin{table}[b!]
\caption{Execution time and Gflop/s of Laplace P2P kernel with SSE and AVX. The calculation is a direct summation between $N=10^5$ particles for both potential and force, which results in 20 floating-point operations per interaction.}
\label{tab:avx}
\begin{center}
\begin{tabular}{|c|c|c|c|}
\hline
& no SIMD & SSE & AVX\\
\hline \hline
CPU time [s] & 1.75 & 0.486 & 0.274\\
\hline
Gflop/s & 48.7 & 189 & 321\\
\hline
Speedup & $\times 1$ & $\times 3.88$ & $\times 6.59$\\
\hline
\end{tabular}
\end{center}
\end{table}

We first demonstrate the performance increase in a Laplace P2P kernel when SSE and AVX intrinsics are used. The runs are on a Intel Xeon E5-2643 (Sandy Bridge, 3.3 GHz, 4 physical cores) with \texttt{g++ -O3 -mavx -fopenmp -ffast-math -funroll-loops}. The execution time and Gflop/s are shown in Table \ref{tab:avx}. The calculation is a direct summation between $N=10^5$ particles for both potential and force, which results in 20 floating-point operations per interaction. In all three cases, the \texttt{\#pragma omp parallel for} directive is used in the outermost loop to utilize all cores. We obtain a speedup of $\times 3.88$ with SSE and $\times 6.59$ with AVX. These SSE and AVX kernels are available from \texttt{https://bitbucket.org/rioyokota/exafmm-dev} as P2P kernels in our FMM. In this code, predefined macros will detect the availability of AVX/SSE and turn them on if the compiler supports them. Furthermore, the loop counter is carried across AVX, SSE, and non-SIMD kernels so that the loop is first processed in increments of 8, then remainders of 8 are processed with SSE, and finally remainders of 4 are processed with non-SIMD kernels.

\subsection{Comparison between different FMM codes.}
\begin{figure}[t]
\centering
\includegraphics[width=0.8\textwidth]{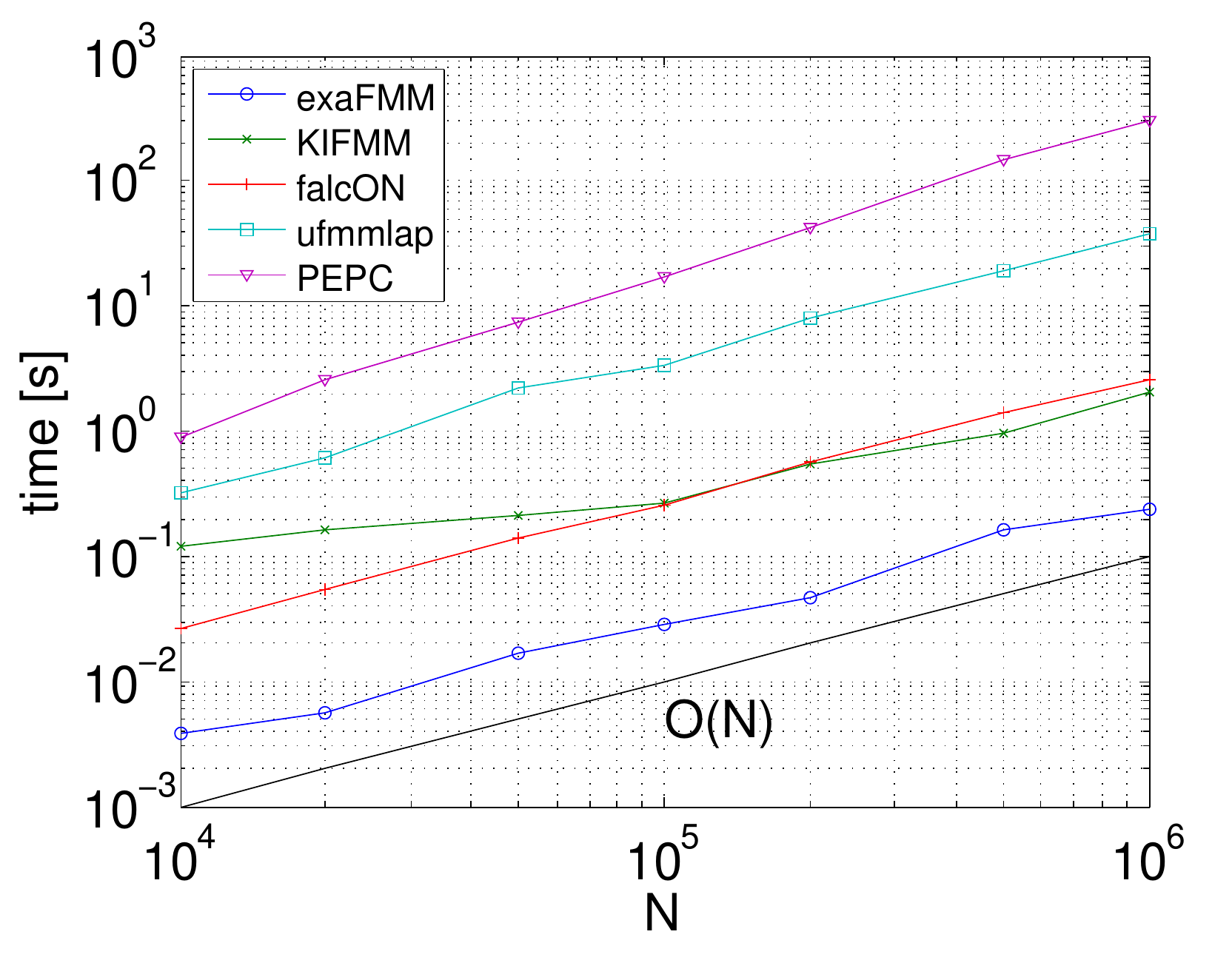}
\caption{Comparison between different FMM codes. The benchmark is for a Laplace kernel with randomly distributed particles in a cube. These benchmarks are performed on a single socket of a Xeon X5650 (Westmere-EP, 2.67GHz, 6 physical cores). Parameters in each code are adjusted to obtain 3 significant digits of accuracy in the force for all codes. Timings include the tree construction time.}
\label{fig:comparison}
\end{figure}

In order to clarify the relative performance of the major open source treecodes and FMMs, we perform a benchmark of 5 different codes on the same computer. The benchmark is for a Laplace kernel with randomly distributed particles in a cube. The number of particles is changed from $N=10^4$ to $10^6$, while setting the accuracy to 3 significant digits for the force (not potential). Note that this setting will yield a little more than 4 significant digits for the potential. For codes that do not calculate the force, we tried our best to match the corresponding accuracy for the potential. These benchmarks are performed on a single socket of a Xeon X5650 (Westmere-EP, 2.67GHz, 6 physical cores). \texttt{exaFMM} is the latest version of our code\footnote{https://bitbucket.org/rioyokota/exafmm-dev}. \texttt{KIFMM} \cite{Chandramowlishwaran2010} is the latest version directly obtained from the developers. \texttt{falcON} \cite{Dehnen2002} is publicly available as part of the package \texttt{NEMO}\footnote{http://bima.astro.umd.edu/nemo}. \texttt{ufmmlap} \cite{Greengard1997} is publicly available as \texttt{FMM-Laplace-U} from \texttt{fastmultipole.org}\footnote{http://fastmultipole.org/Main/FMMSuite}. \texttt{PEPC} \cite{Winkel2012} is the latest publicly available version\footnote{https://trac.version.fz-juelich.de/pepc}.

\begin{table}[b!]
\caption{Characteristics of the different codes in Figure \ref{fig:comparison}. DTT : dual tree traversal, STT : single tree traversal, EC : equivalent charges, TBB : thread building blocks}
\label{tab:characteristics}
\begin{center}
\begin{tabular}{|c|c|c|c|c|c|}
\hline
& \texttt{exaFMM} & \texttt{KIFMM} & \texttt{faclON} & \texttt{ufmmlap} & \texttt{PEPC}\\
\hline \hline
Language & C++ & C++ & C++ & F90 & Fortran 2003\\
\hline
Tree data structure & linked list & array & linked list & array & linked list\\
\hline
Interaction list & DTT & U,V,W,X \cite{Ying2004} & DTT & U,D,N,S,E,W \cite{Greengard1997} & STT\\
\hline
Series expansion & Cartesian & EC + FFT & Cartesian & planewave & Cartesian\\
\hline
Accuracy control & $p$ \& $\theta$ & $p$ & $\theta$ & $p$ & $\theta$\\
\hline
MPI & yes & yes & no & no & yes\\
\hline
Threading model & TBB & OpenMP & N/A & N/A & pthreads\\
\hline
SIMD & AVX & SSE & N/A & N/A & N/A\\
\hline
\end{tabular}
\end{center}
\end{table}

The characteristics of the 5 different codes are summarized in Table \ref{tab:characteristics}. \texttt{PEPC} is a treecode, while the other four are FMMs. \texttt{PEPC} inherits most of it's features from Salmon and Warrens treecode \cite{Warren1993}, and has the structure of a typical treecode. On the other hand, \texttt{ufmmlap} can be thought of as a typical FMM code, so the contrast in the tree data structure, interaction list, series expansion, and accuracy control are most drastic between these two. \texttt{KIFMM} is similar to \texttt{ufmmlap} in the sense that it also uses a typical FMM framework, except the series expansion is unique. The high performance of \texttt{KIFMM} mostly comes from the implementation (OpenMP+SSE with low level optimizations), and on a 12 (hyper-threaded) core machine this results in a big difference between the single core non-SIMD implementation of \texttt{ufmmlap}. Although \texttt{falcON} is also a single core non-SIMD implementation, it achieves the same speed as \texttt{KIFMM} just from purely algorithmic improvements. These algorithmic improvements include the dual tree traversal, mutual interaction, error aware local optimization of the opening angle $\theta$, and a multipole acceptance criterion based on $B_{max}$ and $R_{max}$ \cite{Dehnen2002}. Finally, \texttt{exaFMM} combines the algorithmic improvements of \texttt{falcON} with a highly parallel implementation using MPI, Intel TBB, and AVX. Accuracy can be controlled by both $p$ and $\theta$, the dual tree traversal is parallelized efficiently with task-based parallelism that can either use OpenMP 3.0 or Intel TBB, and the kernels are optimized with AVX to yield 300 Gflop/s on a single CPU socket. With the combination of these techniques, \texttt{exaFMM} can calculate the Laplace potential+force for $N=10^6$ particles to 3 digits of accuracy for the force in approximately $0.2$ seconds on a single CPU socket, which exceeds the performance of most GPU implementations of FMM.

\section{Summary}
The present work attempts to integrate the independent efforts in the fast N-body community to create the fastest N-body library for many-core architectures of the future. Focus is placed on low accuracy optimizations, in response to the recent interest to use FMM as a preconditioner for sparse linear solvers. A direct comparison with other state-of-the-art fast $N$-body codes demonstrates that orders of magnitude increase in performance can be achieved by careful selection of the optimal algorithm and low-level optimization of the code.

The low Byte/flop of future architectures affects the choice of series expansions and translation methods, and our focus on low accuracy perturbs this decision even further. For example, the Cartesian Taylor expansion was shown to be much faster that spherical harmonics with rotation in the range of accuracy of interest. Furthermore, simultaneous control over the order of expansion $p$ and opening angle $\theta$ was also shown to yield a significant increase in performance over conventional treecodes and FMMs that only control either of them.

The dual tree traversal offers an interesting alternative to conventional interaction list construction, and modern compilers now provide the necessary tools to parallelize them efficiently. The fact that dual tree traversals allow rectangular cells, enables great flexibility in the partitioning/load-balancing on distributed memory architectures. There has been little work in this promising area, because dual tree traversal has not been used in the FMM community until now. Although the scope of the present article is focused on single node performance, we would like to point out that the current framework enables various interesting possibilities for partitioning/load-balancing FMMs.

The performance of the FMM code described in this article is made available to the general public\footnote{https://bitbucket.org/rioyokota/exafmm-dev}. The developer has made a significant effort to make the code readable by carefully designing the code structure to make it as simple as possible, while retaining the performance. One can easily see by comparing the source code with the other 4 codes in the benchmark study, that \texttt{exaFMM} is more than an order of magnitude smaller in terms of lines of code, while being orders of magnitude faster.

\end{document}